\documentclass[12pt]{amsart}
\usepackage{amsmath,amsthm,amscd,amsfonts,amssymb}

\usepackage{amsmath,amsthm,amscd,amsfonts,amssymb}
 1   

\newcommand{\bN}{{\mathbb{N}}}
\newcommand{\bZ}{{\mathbb{Z}}}

\newcommand{\cD}{{\mathcal D}_R}
\newcommand{\cmb}[2]{\footnotesize \left (\begin{array}{c} #1 \\ #2 \end{array} \right )}

\newtheorem{theorem}{Theorem}[section]
\newtheorem{lemma}[theorem]{Lema}

\theoremstyle{definition}

\newtheorem{example}[theorem]{Example}

\theoremstyle{remark}
\newtheorem{remark}[theorem]{Remark}

\numberwithin{equation}{section}

\begin{document}

\title[A surprising fact about $\mathcal{D}$-modules in characteristic $p>0$]
{A surprising fact about $\mathcal{D}$-modules in characteristic
$p>0$}

\date{\today}

\author{Josep \`Alvarez Montaner}
\thanks{Research of the first author partially supported by a
Fulbright grant and the Secretar\'{\i}a de Estado de Educaci\'on y
Universidades of Spain and the European Social Funding}

\address{Departament de Matem\`atica Aplicada I\\
Universitat Polit\`ecnica de Catalunya\\ Avinguda Diagonal 647,
Barcelona 08028, SPAIN} \email{Josep.Alvarez@upc.es}

\author{Gennady Lyubeznik}
\thanks{The second author greatfully
acknowledges NSF support}

\address{Department of Mathematics\\
University of Minnesota\\ 206 Church St. S.E., Minneapolis, MN
55455, USA} \email{gennady@math.umn.edu}

\maketitle


\begin{abstract}
Let $R=k[x_1, \dots ,x_d]$ be the polynomial ring in $d$
independent variables, where $k$ is a field of characteristic
$p>0$. Let $\cD$ be the ring of $k$-linear differential operators
of $R$ and let $f$ be a polynomial in $R$. In this work we prove
that the localization $R[\frac{1}{f}]$ obtained from $R$ by
inverting $f$ is generated as a $\cD$-module by $ \frac{1}{f}$.
This is an amazing fact considering that the corresponding
characteristic zero statement is very false.

\end{abstract}

\thispagestyle{empty}

\section{Introduction}

Let $k$ be a field and let $R=k[x_1,\dots, x_d]$, or
$R=k[[x_1,\dots, x_d]]$ be either a ring of polynomials or formal
power series in a finite number of variables over $k$. Let $\cD$
be the ring of $k$-linear differential operators on $R$. For every
$f\in R$, the natural action of $\cD$ on $R$ extends uniquely to
an action on $R[\frac{1}{f}]$ via the standard quotient rule.
Hence $R[\frac{1}{f}]$ acquires a natural structure of
$\cD$-module. It is a remarkable fact that $R[\frac{1}{f}]$ has
finite length in the category of $\cD$-modules. This fact has been
proven in characteristic $0$ by Bernstein \cite[Corollary
1.4]{Be72} in the polynomial case and by Bj\"ork \cite[Theorems
2.7.12, 3.3.2]{Bj79} in the formal power series case and in
characteristic $p>0$ by B$\o$gvad \cite[Proposition 3.2]{Bo95} in
the polynomial case and by Lyubeznik \cite[Theorem 5.9]{Ly97} in
the formal power series case. Thus the ascending chain of
submodules
$$  \cD \cdot \frac{1}{f} \subseteq \cD \cdot \frac{1}{f^2}
\subseteq \cdots \subseteq \cD \cdot \frac{1}{f^k} \subseteq
\cdots \subseteq R[\frac{1}{f}]  $$ stabilizes, i.e.
$R[\frac{1}{f}]$ is generated by $\frac{1}{f^i}$ for some $i$.
This paper is motivated by the following natural question: {\it
What is the smallest $i$ such that $\frac{1}{f^i}$ generates
$R[\frac{1}{f}]$ as a $\cD$-module?}

When $k$ is a field of characteristic zero and $f\in R$ is a
non-zero polynomial it has been proven in \cite[Theorem 1']{Be72}
that there exists of a monic polynomial $b_f(s)\in k[s]$ and a
differential operator $Q(s) \in \cD[s]$ such that $$ Q(s)\cdot
f^{s+1}=b_f(s)\cdot f^s$$ for every $s$. The polynomial $b_f(s)$
is called the Bernstein-Sato polynomial of $f$ and is always a
multiple of $(s+1)$.  Let $-i$ be the least integer root of
$b_f(s)$. Then, $b_f(s)\ne 0$ for any integer $s<-i$ and therefore
$f^s\in \cD\cdot f^{s+1}$. In particular, $R[\frac{1}{f}]$ is
generated by $\frac{1}{f^i}$ and, as is shown in \cite[Lemma
1.3]{Wa03}, it cannot be generated by $\frac{1}{f^j}$ for $j<i$.
This gives a complete answer to our question in characteristic
zero.

For example, consider the polynomial $f=x_1^2+\cdots+x_{2n}^2 \in
R=k[x_1,\dots,x_{2n}]$. Then we have the functional equation
$$\frac{1}{4}(\frac{\partial^2}{\partial x_1^2} +
\cdots + \frac{\partial^2}{\partial x_{2n}^2} )\cdot f^{s+1} =
(s+1)(s+n) \cdot f^s$$ where the Bernstein-Sato polynomial is
$b_f(s)=(s+1)(s+n)$. Hence $R[ \frac{1}{f}]$ is generated by $
\frac{1}{f^n}$ as a $\cD$-module but it can not be generated by $
\frac{1}{f^i}$ for $i<n$.

But in characteristic $p>0$ a differential operator of a fixed
order annihilates the powers $\frac{1}{f^{p^{s}}}$ for $s$ large
enough, so a functional equation such as above even if it exists
does not imply that $\frac{1}{f^{p^{s}+1}}\in \cD \cdot
\frac{1}{f^{p^s}}$ for all $s$.

The goal of this paper is to prove the following amazing result.

\begin{theorem}\label{T1}
Let $R=k[x_1,\dots,x_d]$ where $k$ is a field of characteristic
$p>0$ and let $f\in R$ be a non-zero polynomial. Then
$R[\frac{1}{f}]$ is generated by $\frac{1}{f}$ as $\cD$-module.
\end{theorem}

Our proof does not extend to the case of formal power series; some
new idea seems to be needed in this case (see Remark \ref{R1}).

\noindent {\bf Acknowledgement} The first author would like to
thank the Department of Mathematics at the University of Minnesota
for the warm welcome he received during his postdoctoral stay.

\section{Differential operators in positive characteristic}


Let $\bN$ be the set of non-negative integers. Throughout, we will
use multi-index notation in the polynomial ring $R=k[x_1, \dots
,x_d]$, where $k$ is a field of characteristic $p>0$. So, given
$\alpha=(\alpha_1,\dots,\alpha_d)\in \bN^d$ we will denote the sum
of its components by $|\alpha|$ and ${\underline{x}^{\alpha}}$
will stand for the monomial
${\underline{x}^{\alpha}}=x_1^{\alpha_1}\cdots x_d^{\alpha_d}$. A
pair of multi-indices $\alpha$ and $\beta$ are ordered as usual:
$\alpha < \beta$ if and only if $\alpha_i < \beta_i$ for
$i=1,\dots,d$.

For a general description of the ring of differential operators we
refer to \cite[$\S$16.8]{GD67}. For the case we are considering in
this work we refer to \cite[Th\'eor\`eme 16.11.2]{GD67}. The ring
of differential operators $\cD=D(R,k)$ associated to the
polynomial ring $R$ is the ring extension of $R$ generated by the
set of differential operators
$$\{\hskip 2mm D_{t,i}=\frac{1}{t!}\frac{\partial^{t}}{\partial x_i^{t}}
\hskip 2mm | \hskip 2mm t\in \bN \hskip 2mm,\hskip 2mm i=1,\dots,d
\hskip 2mm \}$$

Given $\beta \in \bN^d$, ${D_{\beta}}$ will denote the
differential operator ${D_{\beta}}:= D_{\beta_1,1}\cdots
D_{\beta_d,d}$. We can extend the multi-index notation to $\cD$
considering the $k$-basis formed by the monomials
$\underline{x}^\alpha D_\beta$. In the sequel, a differential
operator $Q\in \cD$ will be written in right normal form, i.e.
$$Q =\sum_{\alpha,\beta\in \bN^d} a_{\alpha \beta}
\hskip 1mm {\underline{x}^\alpha} \hskip 1mm D_\beta,$$ where all
but finitely many $a_{\alpha \beta}\in k$ are zero.

Let $k[R^{p^n}]$ be the $k$-algebra generated by $p^n$-th powers
of elements of $R$. Let $\cD^{(n)}$ be the ring extension of $R$
generated by the set of differential operators up to order $p^n$,
i.e. $\{\hskip 2mm D_\alpha \hskip 2mm | \hskip 2mm \alpha <
{\underline{p}^n} \hskip 2mm \} $ where
${\underline{p}^n}=(p^n,\dots,p^n)\in \bN^d$. Then we have an
increasing chain of finitely generated ring extensions of $R$
$$\cD^{(0)} \subseteq \cD^{(1)} \subseteq \cD^{(2)} \subseteq \cdots
\subseteq \cD$$ whose union is $\cD$.

\begin{lemma}\label{L2}
Let $Q\in \cD^{(n)}$ and $f\in k[R^{p^n}]$. Then $Q$ commutes with
$f$, i.e.  $Q(f\cdot g)= f \cdot Q(g)$ for all $g\in R$.
\end{lemma}

\begin{proof}
Writing out $Q$, $f$ and $g$ as sums of monomials and considering
that $D_{t,i}$ commutes with $D_{t',j}$ and $x^s_j$ for $j\neq i$,
one sees that it is enough to prove the statement for $Q=D_{t,i}$,
$f=x_i^{p^n}$ and $g=x_i^s$, where $t<p^n$, $s$ is an integer and
$i=1,\dots,d$. In this case we have
$$D_{t,i}(x_i^{p^n+s})= x_i^{p^n} \cdot D_{t,i}(x_i^s)$$ just comparing
the coefficient at $x_i^{p^n+s-t}$ on both sides modulo $p$.
\end{proof}


\section{Proof of Theorem \ref{T1}}

We notice first that it is enough to show that if $f\in R$ is a
non-zero polynomial, then $\frac{1}{f^p}$ belongs to the
$\cD$-module generated by $\frac{1}{f}$. Once this is proved we
can apply this result to $f^{p^{s-1}}$ to get
$\frac{1}{f^{p^s}}\in \cD\cdot \frac{1}{f^{p^{s-1}}}$ for every
$s>1$. Since the set $\frac{1}{f}, \frac{1}{f^p},
\frac{1}{f^{p^2}}, \dots $ generates $R[\frac{1}{f}]$ as
$R$-module, we are done.

We can also reduce to the case of $k$ being a perfect field. If
$k$ is not perfect, let $K$ be the perfect closure of $k$. Assume
there is a differential operator $Q=\sum a_{\alpha \beta}
{\underline{x}^\alpha} D_\beta$ with coefficients $a_{\alpha
\beta}\in K$ such that $Q(\frac{1}{f})= \frac{1}{f^p}$. This is
equivalent to the fact that a system of finitely many linear
equations with coefficients in $k$ has solutions in $K$, where the
non-zero coefficients $a_{\alpha\beta}$ of $Q$ are thought of as
the unknowns of the system. For example, if $f=x_1$, we may be
looking for a solution in the form $Q=a D_{p-1,1}$, so we get an
equation $Q(\frac{1}{x})=\frac{1}{x^p}$. Since $Q(\frac{1}{x})=a
\frac{1}{x^p}$, the corresponding linear system is just one
equation $a=1$. The system has a solution in $K$, namely, the
coefficients of $Q$. Hence it is consistent, so it must have a
solution in $k$ because the coefficients of the linear system are
in $k$ (in fact the coefficients are in the prime field $\bZ /
p\bZ$). So there is a differential operator $Q'$ with coefficients
in $k$ such that $Q'(\frac{1}{f})= \frac{1}{f^p}$.

Henceforth we will assume that the base field $k$ of our
polynomial ring is perfect. It is enough to show that under this
assumption $ \frac{1}{f^p}$ belongs to the $\cD$-submodule
generated by $ \frac{1}{f}$, for all $f\in R$.

Given a polynomial $f\in R$ and an integer $n\geq 1$, we can write
in a unique way
$$f({\underline{x}})= \sum_{0\leq \alpha < {\underline{p}^n}}
f_{\alpha}({\underline{x}^{\underline{p}^n}})\hskip 2mm
{\underline{x}^\alpha}$$ where $f_{\alpha}({\underline{z}}) \in
k[{\underline{z}}]$ are polynomials in $d$ variables.  Consider
the ideal $J_n(f)$ generated by the polynomials
$f_{\alpha}(\underline{x}^{\underline{p}^n})$ in the decomposition
of $f$ with respect to ${\underline{x}}^{\underline{p}^n}$.

\begin{lemma} \label{L4}
Let $f,g,h\in R$ be polynomials such that $f=g\cdot h$. Then
$J_n(f)\subseteq J_n(g)$.
\end{lemma}

\begin{proof}
Consider the decomposition of $g$ with respect to
${\underline{x}}^{\underline{p}^n}$
$$g({\underline{x}})= \sum_{0\leq \alpha < {\underline{p}^n}}
g_{\alpha}({\underline{x}^{\underline{p}^n}})\hskip 2mm
{\underline{x}^\alpha}$$ Set $h({\underline{x}})=\sum a_\beta
\hskip 2mm \underline{x}^\beta$, then $$f({\underline{x}})=
g({\underline{x}})\cdot h({\underline{x}})= \sum_{0\leq \alpha <
{\underline{p}^n}}
g_{\alpha}({\underline{x}^{\underline{p}^n}})\hskip 2mm (\sum
a_\beta \hskip 2mm{\underline{x}^{\alpha+\beta}})$$ Rewriting  in
the form
$$f({\underline{x}})= \sum_{0\leq \gamma < {\underline{p}^n}}
f_{\gamma}({\underline{x}^{\underline{p}^n}})\hskip 2mm
{\underline{x}^\gamma}$$ we get
$$f_{\gamma}({\underline{x}^{\underline{p}^n}})=\sum a_\beta \hskip 2mm g_{\alpha}({\underline{x}^{p^n}})
\hskip 2mm \underline{x}^{j \underline{p}^n}$$ where the sum is
taken over the multi-indices $\alpha$ and $\beta$ such that
$\underline{x}^{\alpha+\beta}= \underline{x}^{j \underline{p}^n}
\underline{x}^{\gamma}$ for a given $j\in \bN$. In particular
$f_{\gamma}({\underline{x}^{\underline{p}^n}})\in J_n(g)$.
\end{proof}

\begin{lemma} \label{L5}
Let $f,g\in R$ be polynomials such that $f=g^p$. Then $J_n(f)=
J_{n-1}(g)^{[p]}$

\end{lemma}

\begin{proof}
It is enough to raise to the $p-$th power the decomposition of $g$
with respect to ${\underline{x}}^{\underline{p}^{n-1}}$.
\end{proof}


Notice that $\cD^{(n)} \cdot f$ is an ideal of $R$. This ideal can
be also described as follows:

\begin{lemma} \label{L6}
 $J_n(f)=\cD^{(n)} \cdot f$.
\end{lemma}

\begin{proof}
By Lemma \ref{L2} every $Q \in \cD^{(n)}$ commutes with every
$f_{\alpha}(\underline{x}^{\underline{p}^n})$ in the decomposition
of $f$ with respect to ${\underline{x}}^{\underline{p}^n}$. Hence
$$Q(f)=\sum_{0\leq \alpha < {\underline{p}^n}}
f_{\alpha}({\underline{x}^{\underline{p}^n}})\hskip 2mm
Q({\underline{x}^\alpha})$$ In particular, $Q(f)\in J_n(f)$, i.e.
$\cD^{(n)} \cdot f \subseteq J_n(f)$. To prove the opposite
containment it is enough to show that every
$f_{\alpha}(\underline{x}^{\underline{p}^n})$ belongs to the ideal
$\cD^{(n)} \cdot f$.

Consider the multi-index $\beta=(p^n-1,\dots,p^n-1)\in \bN^d$.
Then we have
$$D_\beta (f)= f_\beta(\underline{x}^{\underline{p}^n})\in \cD^{(n)} \cdot f$$
Now we proceed by induction on
$\sigma_\beta=\Sigma_{i=1}^{d}(p^n-1-\beta_i)$, the case
$\sigma_\beta=0$ being just proved. Let $\beta\in \bN^d$ be a
multi-index such that $D_\beta \in \cD^{(n)}$. Then we have
$$D_\beta (f)= f_\beta(\underline{x}^{\underline{p}^n}) +
\sum_{\beta' > \beta } f_{\beta'}(\underline{x}^{\underline{p}^n})
\hskip 1mm (a_{\beta'\beta}\hskip 1mm
\underline{x}^{\beta'-\beta})$$ where $a_{\beta'\beta}=
\cmb{\beta_1'}{\beta_1} \cdots \cmb{\beta'_d}{\beta_d}$. By
induction $f_{\beta'}(\underline{x}^{\underline{p}^n})\in
\cD^{(n)}\cdot f$ for any $\beta'>\beta $, so we are done.
\end{proof}

Since $k$ is a perfect field, the coefficients of
$f_{\alpha}({\underline{x}^{\underline{p}^n}})$ in the
decomposition of $f$ with respect to
${\underline{x}}^{\underline{p}^n}$ are $p^n$-th powers, hence
$f_{\alpha}({\underline{x}^{\underline{p}^n}})=(\widetilde{f}_{\alpha}({\underline{x}}))^{p^n}$,
where $\widetilde{f}_{\alpha}({\underline{x}})$ are polynomials in
$R$. Consider the ideal $I_n(f)$ generated by the polynomials
$\widetilde{f}_{\alpha}(\underline{x})$. Notice that $J_n(f)$ is
the $n$-th Frobenius powers of the ideal $I_n(f)$, i.e. $J_n(f)=
I_n(f)^{[p^n]}$.

\begin{lemma}\label{L7}
Let $f\in R$ be a polynomial. For any integer $n\geq 1$ there is
an inclusion of ideals
$$I_n(f^{p^n-1})\subseteq I_{n-1}(f^{p^{n-1}-1})$$

\end{lemma}

\begin{proof}
 The equality $J_n(f^{p^n-p})= J_{n-1}(f^{p^{n-1}-1})^{[p]}$
given by Lemma \ref{L5} translates to
$$I_n(f^{p^n-p})^{[p^n]}=
(I_{n-1}(f^{p^{n-1}-1})^{[p^{n-1}]})^{[p]}=I_{n-1}(f^{p^{n-1}-1})^{[p^n]}$$
This implies $$I_n(f^{p^n-p})= I_{n-1}(f^{p^{n-1}-1})$$ due to the
fact that the polynomial ring $R$ is regular. On the other hand,
the inclusion $J_{n}(f^{p^{n}-1})=J_n(f^{p^{n}-p} \cdot
f^{p-1})\subseteq J_{n}(f^{p^{n}-p})$ given by Lemma \ref{L4}
implies
$$I_{n}(f^{p^{n}-1})\subseteq I_{n}(f^{p^{n}-p})$$ again because $R$ is regular so we
get the desired inclusion.
\end{proof}

\begin{lemma}\label{L8}
The descending chain of ideals
$$I_1(f^{p-1})\supseteq I_2(f^{p^2-1})\supseteq \cdots \supseteq I_{n-1}(f^{p^{n-1}-1})\supseteq I_n(f^{p^n-1})\supseteq \cdots$$
stabilizes.
\end{lemma}

\begin{proof}
Assume that $\deg f= e$. Let
$$f^{p^n-1}({\underline{x}})= \sum_{0\leq \alpha <
{\underline{p}^n}}
f_{\alpha}({\underline{x}^{\underline{p}^n}})\hskip 2mm
{\underline{x}^\alpha}$$ be the decomposition of $f^{p^n-1}$ with
respect to ${\underline{x}}^{\underline{p}^n}$. Since $\deg
f^{p^n-1}=e(p^n-1)$, the polynomials in the decomposition satisfy
$\deg f_{\alpha}({\underline{x}^{\underline{p}^n}})\leq e(p^n-1)<
ep^n$. On the other hand, $\deg
f_{\alpha}({\underline{x}^{\underline{p}^n}})= p^n \deg
\widetilde{f}_{\alpha}({\underline{x}})$ implies $\deg
\widetilde{f}_{\alpha}({\underline{x}})< e$.

Let $W_e$ be the $k$-vector space of polynomials of degree
strictly smaller than $e$. For every $n$, the ideal
$I_n(f^{p^n-1})$ is generated by $I_n(f^{p^n-1})\cap W_e$. Thus we
have a descending sequence of $k$-vector subspaces of $W_e$
$$(I_1(f^{p-1})\cap W_e)\supseteq (I_2(f^{p^2-1})\cap
W_e)\supseteq \cdots \supseteq (I_n(f^{p^n-1})\cap W_e)\supseteq
\cdots$$ that must stabilize because $W_e$ is a finite-dimensional
$k$-vector space.
\end{proof}

Now we can complete the proof of Theorem \ref{T1} as follows.
Assume that the chain of ideals given in Lemma \ref{L8} stabilize
at the level $s-1$, i.e. $I:=I_{s-1}(f^{p^{s-1}-1})=
I_s(f^{p^s-1})$. From the equalities $J_s(f^{p^s-1})=I^{[p^s]}$
and $J_{s-1}(f^{p^{s-1}-1})=I^{[p^{s-1}]}$ we deduce
$$ J_s(f^{p^s-p})= J_{s-1}(f^{p^{s-1}-1})^{[p]}= (I^{[p^{s-1}]})^{[p]}=I^{[p^{s}]}=J_s(f^{p^s-1})$$
Thus we have $f^{p^s-p}\in J_s(f^{p^s-1})$. By Lemma \ref{L6}
there is a differential operator $Q\in \cD^{(s)}$ such that
$Q(f^{p^s-1})=f^{p^s-p}$. Since $Q$ commutes with $f^{p^s}$, we
see that
$$Q(\frac{f^{p^s-1}}{f^{p^s}})= \frac{f^{p^s-p}}{f^{p^s}}$$ so we
get $Q( \frac{1}{f})=\frac{1}{f^p}$ as we desired. \hskip 7cm
$\Box$

\begin{remark}\label{R1}
In the case of formal power series all parts of our proof go
through except the proof of Lemma \ref{L8}. Most likely, the
statement of Lemma \ref{L8} is still true in the case of formal
power series but a very different proof is needed.
\end{remark}

\begin{example}
Let $R=k[x_1,x_2,x_3,x_4]$ where $k$ is a field of characteristic
$p>0$. Consider the polynomial $f=x_1^2+x_2^2+x_3^2+x_4^2 $. We
are going to find a differential operator $Q\in \cD$ such that $Q(
\frac{1}{f})=\frac{1}{f^p}$ just checking out the monomials in
$f^{p-1}$. Let $S(f,p,n)$ be the set of terms $a_\alpha
\underline{x}^\alpha$ in $f^{p-1}$ such that $\alpha <
\underline{p}^n$.

The set $S(f,p,1)$ is non-empty. Namely we have:

$\bullet$ If $4$ divides $p-1$, then $a_\alpha
\underline{x}^\alpha \in S(f,p,1)$ where

\hskip 1cm
$\alpha=(\frac{p-1}{2},\frac{p-1}{2},\frac{p-1}{2},\frac{p-1}{2})$

\hskip 1cm $a_\alpha= \frac{(p-1)!}{(\frac{p-1}{4}!)^4}$

$\bullet$ If $4$ does not divide $p-1$, then $a_\alpha
\underline{x}^\alpha \in S(f,p,1)$ where

\hskip 1cm
$\alpha=(\frac{p+1}{2},\frac{p+1}{2},\frac{p-3}{2},\frac{p-3}{2})$

\hskip 1cm $a_\alpha=
\frac{(p-1)!}{(\frac{p+1}{4}!)^2(\frac{p-3}{4}!)^2}$

\noindent Let $a_\alpha \underline{x}^\alpha $  be a leading term
of $S(f,p,1)$ with respect to the usual order. Notice that
$\frac{1}{a_\alpha} \hskip 1mm D_\alpha(f^{p-1})=1$. The
differential operator $Q= \frac{1}{a_\alpha} \hskip 1mm D_\alpha$
commutes with $f^p$ by Lemma \ref{L2} so we get the desired
result.

\end{example}


\nocite{*}


\begin{thebibliography}{22}

\bibitem[1]{Be72}
I.~N.~Bern$\check{s}$te$\breve{\i}$n,
\newblock  Analytic continuation of generalized functions with respect to a parameter,
\newblock Funkcional. Anal. i Prilo$\check{z}$en., 6 (4) (1972) 26--40.


\bibitem[2]{Bj79}
J.~E. Bj{\"o}rk,
\newblock Rings of differential operators, North
  Holland Mathematics Library,
  \newblock Amsterdam, 1979.



\bibitem[3]{Bo95}
R.~B$\o$gvad,
\newblock  Some results on $D$-modules on Borel varieties in characteristic $p>0$,
\newblock  J. Algebra,  173 (3) (1995) 638--667.




\bibitem[4]{GD67}

A.~Grothendieck and J.~Dieudonn\'e,
\newblock \'El\'ements de g\'eom\'etrie alg\'ebrique IV. \'Etude locale des sch\'emas et des morphismes
de sch\'emas,
\newblock Publications Math\'ematiques I.H.E.S.
32 (1967).




\bibitem[5]{Ly97}
G.~Lyubeznik,
\newblock {$F$-modules: applications to local cohomology and
$D$-modules in characteristic $p>0$,}
\newblock { J. Reine Angew. Math.} { 491} (1997), 65--130.

\bibitem[6]{Wa03}
U.~Walther,
\newblock { Bernstein-Sato polynomials versus cohomology of the Milnor fiber
       for generic hyperplane arrangements,}
\newblock to appear in Compositio Math.




\end{thebibliography}
\end{document}